\documentclass{amsart}

\usepackage{graphicx}

\usepackage{amsmath}
\usepackage{amssymb}
\usepackage{amsbsy}

\newcommand{\trp}{^{\mbox{\scriptsize\sf T}}}
\newcommand{\reals}{{\mathbb R}}

\newtheorem{theorem}{Theorem}[section]

\newtheorem{proposition}[theorem]{Proposition}
\theoremstyle{definition}
\newtheorem{definition}[theorem]{Definition}

\theoremstyle{remark}

\numberwithin{equation}{section}

\begin{document}

\title{A network model for granular statics with impenetrability constraints}


\author{K. A. Ariyawansa}
\address{Department of Mathematics, Washington State University, Pullman, WA, 99164}
\email{ari@math.wsu.edu}
\thanks{Work of K. A. Ariyawansa was supported in part by ARO grant DAAD 19-00-1-0465}

\author{Leonid Berlyand}
\address{Department of Mathematics and Materials Research Institute , Penn State
University, University Park, PA 16802, USA}
\email{berlyand@math.psu.edu}
\thanks{Work of Leonid Berlyand was supported in part by NSF grant DMS-0204637. }

\author{Alexander Panchenko}
\address{Department of Mathematics, Washington State University, Pullman, WA, 99164}
\email{panchenko@math.wsu.edu}
\thanks{Work of Alexander Panchenko was supported in part by DOE grant DE-FG02-05ER25709.}


\date{May 31, 2006}

\begin{abstract}
We study quasi-static deformation of dense granular packings.
In the reference configuration, a granular material is under confining
stress (pre-stress). Then the packing is deformed by imposing external boundary conditions,
which model engineering experiments such as shear and compression. The deformation
is assumed to preserve the local structure of neighbors for each particle, which is a
realistic assumption for highly compacted packings driven by small boundary displacements. 
We propose a two-dimensional network model of such deformations.
The model takes into account elastic interparticle interactions and incorporates
geometric impenetrability constraints. The effects of friction are neglected. 
In our model, a granular packing is represented by a spring-lattice network, whereby
the particle centers correspond to vertices of the network, and interparticle contacts
correspond to the edges. We work with general network geometries: periodicity is not
assumed. For the springs, we use a quadratic elastic energy function. Combined with the linearized impenetrability constraints, this function provides a regularization of the hard-sphere potential for small displacements.

When the network deforms, each spring either preserves its length  (this corresponds
to a solid-like contact), or expands (this represents a broken contact). Our goal is to study 
distribution of solid-like contacts in the energy-minimizing configuration.
We prove that under certain geometric conditions on the network, there are at least two 
non-stretched springs attached to each node, which means that every particle has at least
two solid-like contacts. The result implies that a particle cannot loose contact with all of its neighbors.
This eliminates micro-avalanches as a mechanism for structural weakening in small shear
deformation.

\end{abstract}

\maketitle



\section{Introduction}
\label{Sect:intro}

Materials that are composed of collections of separate, macroscopic solid grains 
belong to the general classification of {\it granular materials}. Examples of such materials are common, including sand, gravel, medicinal pills, coins, and breakfast cereal.
Granular media are important to numerous industries ranging 
from mining to pharmaceuticals. In geophysics, granular 
materials are a central problem in understanding
the physics of earthquakes and tectonic faulting. 
Earthquake fault zones produce granular wear material continuously as a
function of shear and grinding between the fault surfaces.  The wear 
material, known as fault gouge, varies in thickness from 10's of cm to 1000 m
and plays a critical role in determining the fault zone frictional
strength, the stability of fault slip, and the size of the rupture 
nucleation dimension.

Granular media display a variety of complex 
static and  dynamic properties that distinguish them from conventional
solids and liquids. 
The complexity of granular media lies primarily in the collective 
properties of a macroscopic number of grains and how they interact with 
each other. The conditions under which a granular medium is stable 
or flows and the nature of this flow depend critically on the 
distributions of grain size and shape as well as the interactions 
between the grains. The practical importance of granular media combined
with the richness of their physical properties has led to a great deal
of interest from both theoretical and experimental points of view
\cite{DEG, jaeger, kadanoff}.

An important class of granular materials consists of nearly rigid particles that possess the following property: if a moderate force is applied, the particles start to move, and only after a substantial increase of the force the particles deform significantly. In other words, for loads which are not very high, the deformations inside a particle
are small compared with the displacement of the particle center of mass. 
Consequently, the particle shapes change very little, so that
each particle can be associated with a region of space that is inaccessible
to any other particle. This gives rise to constraints on the admissible positions of particles. These {\it impenetrability constraints} are also known as
geometric, kinematic, and excluded volume constraints.

A physical phenomenon related to appearance of constraints is jamming. A particle is jammed when its motion is completely obstructed by the neighbors, so the whole cluster of neighboring particles can only 
move together as a rigid body. 
The corresponding mathematical notion of rigidity (\cite{W}) can be applied to various physical (sphere packings, frameworks (trusses)), as well as mathematical objects. In particular, an important mathematical object associated with any particle packing is a {\it contact graph} defined as follows:
vertices of this graph are particle centers of mass, while edges represent interparticle contacts.



The simplest physical model that exhibits jamming is a classical hard sphere packing. The particles
in this model are represented by rigid spheres,  and the only interparticle forces are reactions of constraints.  Rigidity of hard sphere packings is studied in \cite{Co}. This problem can be formulated as a problem of detecting rigidity of the {\it cable framework} associated with the contact graph of the packing.  The framework is obtained by replacing edges of the contact graph with the cables, and vertices with flexible hinges. The lengths of the cables can increase but not decrease, which models the impenetrability constraints. Recently, a linear programming algorithm for detecting rigidity in hard sphere packings (equivalently, cable frameworks) was proposed in \cite{D}.  

In this work, we also use {\it bar frameworks}. A bar framework
is obtained from a graph by replacing the edges with rigid bars, and vertices with hinges.
A bar framework and the associated graph are called rigid if the only possible vertex motions correspond to rigid body motions of the whole framework.  We note that both bar and cable frameworks can be associated
to the same graph. To generate the bar framework, the edges of a graph are replaced with rigid bars that can only translate and rotate. In the case of the cable framework, one replaces edges with cables that can either move as rigid bodies or stretch. Thus every motion of a bar framework is also permitted by the cable framework, but the converse is not true in general.  Therefore, it is possible that a bar framework associated to a graph is rigid, while the cable framework corresponding to the same graph is not. Both
cable frameworks and bar frameworks are special case of the so-called tensegrity frameworks studied in \cite{CoW}. In a tensegrity framework, properties of edges can vary, e.g. some edges may be bars, others may be cables or struts (that can shrink, but not stretch).

It appears that the currently available mathematical results \cite{Co, CoW, D} on statics of discrete particle systems with geometric constraints deal only with hard particle packings. 
To the best of authors' knowledge, there are
no results on frictional packings, and even elastic frictionless packings have yet not been studied. 
The present work differs from \cite{Co, CoW, D} in several respects. First, all these studies
deal with rigid particles. We consider a somewhat more realistic situation of geometrically
constrained particles with elastic interactions defined by a quadratic potential energy.
Second, while \cite{Co, CoW, D} focus on jamming, we are interested in generic contact patterns of the the energy minimizing configurations.
The packings that we study are not jammed. Their contact graphs are such that the associated bar framework is rigid, but the packing can still deform when external boundary conditions are applied.

The third difference is in the type of the boundary conditions. The conditions in
\cite{D} are periodic or hard wall conditions. The periodic conditions are commonly used to
minimize influence of the boundaries in the problem. However, presence of walls
is a major factor that determines bulk behavior of granular materials. Therefore it seems better to use boundary conditions corresponding to engineering  and physical experiments, where the walls are rigid and may be moving.

A frequently observed property of granular materials is concentration of the bulk
deformation in thin layers called shear bands. Within a band, the contact forces are
weak, and the relative displacements can be on the order of particle size or larger.
For quasi-static flows driven by small shear rates, 
the corresponding patterns are called micro-bands (\cite{Kuhn1}). In that paper, shear
band structures were studied by means of numerical simulations. The simulations
in \cite{Kuhn1} show that the typical size and number of bands in quasi-static shear depend on the imposed shear rate. For small shear rates, the bands have length and width comparable with the particle size. The distribution of these micro-bands within the material is rather uniform. As the shear rate increases, the band structure exhibits coarsening: the number of bands becomes smaller, and the length of each band increases. For sufficiently large shear rates, a single macroscopic shear band appears.

Here, we are interested in the case of small external boundary conditions. In that case, a micro-band
can be formed by weakening of a single contact, or a small group of neighboring contacts. All weak
contacts form a subnetwork of the whole contact network. Such networks of weak contacts
(corresponding to the micro-band patterns in \cite{Kuhn1}) were studied numerically in \cite{Radjai}. 
The goal of this paper is to describe some generic geometric features of micro-band, or weak contact
networks in dense packings of nearly rigid particles. In small deformation, pattern formation may be caused
by the local jamming (which mathematically amounts to impenetrability constraints), and friction. 
We are concerned with the role of constraints, while friction in neglected. 
The notion of high density at this point is rather intuitive. It could mean, for instance, that each particle is
in contact with at least three other particles, and the packing is jammed in the reference configuration, subject to zero boundary conditions. Below we make the notion of high density more precise (see the second paragraph on p. 4), using the relationship between the contact graph and the Delaunay graph (see e.g. \cite{Ed}), generated by the set of particle centers.

In two dimensions, particles are represented by disks $D_i$ of radii $a_i$ with centers ${\bf x}^i$, $i=1, 2,\ldots, N$. The initial reference configuration is deformed by applying prescribed small displacements to the boundary particles. Assuming that the deformations inside of the individual particles are small, and neglecting rotational degrees of freedom, one can characterize the deformations of  $D_i$ by the displacements ${\bf u}^i$ of their centers. The elastic interaction forces are modeled as
in classical mechanics of point particles: the force exerted by $D_j$ on $D_i$ is applied
at ${\bf x}^i$, its direction is along the line joining ${\bf x}^i$ and ${\bf x}^j$, and its magnitude
depends linearly on ${\bf u}^j-{\bf u}^i$.

We further assume that
the granular material is pre-stressed (or, equivalently, the material is under confining stress). This means that in the reference state, the particles are squashed into each other as a result of applied external pressure.  
Further compression is supposed to be impossible (requires infinite energy), which
introduces impenetrability constraints into the problem.
To model impenetrability, one can, for instance, require that 
\begin{equation}
\label{0.1}
|({\bf x}^i+{\bf u}^i)-({\bf x}^j+{\bf u}^j)|\geq a_i+a_j,
\end{equation}
for each pair of particles. Since the the packing is dense, and the displacements are 
expected to be small, it makes sense to require that a particle cannot escape a cage formed by its neighbors. Therefore, the contacts that exist in the reference configuration may be broken, but no new contacts are created after applying
external  boundary conditions to the reference configuration. 
An important consequence of this assumption is as follows. If the displacements satisfy ({0.1}) for each pair of particles {\it in contact}, then for any pair of particles ({0.1}) is automatically satisfied. Indeed, if two particles are not in contact
in the reference configuration, they cannot come into contact in the deformed configuration, and the distance between
them must be larger than the sum of their radii. In the sequel, we use this assumption in the course of proving the main result of this paper (Theorem \ref{main-theorem}).

Next, we introduce a network model which describes a granular material under the above assumptions. The vertices of the network are the particle centers, and the edges represent
particle contacts. 
The collection of vertices ${\bf x}^i$, $i=1, 2, \ldots, N$, and edges forms the contact network (graph) $\Gamma$. We suppose that $\Gamma$  is a triangulation of a connected, convex polygonal domain $\Omega$. 
This assumption is realistic, since, for example, a periodic 2D packing of disks is triangular. 
Another natural triangulation generated by ${\bf x}^i, i=1, 2,\ldots, N$ is the Delaunay graph $G$. In principle, 
$\Gamma$ and $G$ may be different, since some edges in $G$ may not correspond to contacts. In the present case, we suppose that $\Gamma$ and $G$ coincide, which corresponds to "maximally dense" packings. 

For small displacements ${\bf u}^i$, the quadratic constraints (\ref{0.1}) can be approximated by
their linearizations near ${\bf u}^i={\bf u}^j=0$, which leads to the linearized impenetrability constraints
\begin{equation}
\label{i3-intro}
({\bf u}^j-{\bf u}^i)\cdot {\bf q}^{ij}\geq 0,
~~i=1,2,\ldots, N
\end{equation}
for each pair of vertices $i,j$ connected by an edge of $\Gamma$.
In (\ref{i3-intro}), 
$
{\bf q}^{ij}=({\bf x}^j-{\bf x}^i)/\left|{\bf x}^j-{\bf x}^i\right|
$
are unit vectors
that point from ${\bf x}^i$ to ${\bf x}^j$ along the line of centers.
Note that if the position of $D_i$ is fixed (${\bf u}^i=0$), then
${\bf u}^j$ satisfying (\ref{i3-intro}) must lie in the half-plane ${\bf u}\cdot {\bf q}^{ij}\geq 0$, so that $D_j$ would be moving away from $D_i$. 

For certain boundary conditions the deformed packing can become more loose than the reference packing. On the macroscale, this can be observed as swelling of the specimen caused by the increase in the volume of the
void space between the particles. Such swelling is typical in shear deformation, where the overall volume increase, known as dilatation, is observed in experiments.  To increase the void volume, some of the contacts
present in the reference configuration must be broken in the deformed configuration. Therefore, 
among all the contacts (satisfying (\ref{i3-intro})), we further distinguish two types of contacts: broken and solid-like (see Fig. 1). We call a contact broken if
\begin{equation}
\label{broken}
({\bf u}^j-{\bf u}^i)\cdot{\bf q}^{ij}>0,
\end{equation}
and solid-like if
\begin{equation}
\label{solid-intro}
({\bf u}^j-{\bf u}^i)\cdot{\bf q}^{ij}=0.
\end{equation}

The solid-like contacts correspond to two possible types of pair motions. The first type 
is a rigid motion of a pair, in which case the contact is called {\it stuck}.
\begin{figure}[htbp]
\label{contacts}
\begin{center}
\input{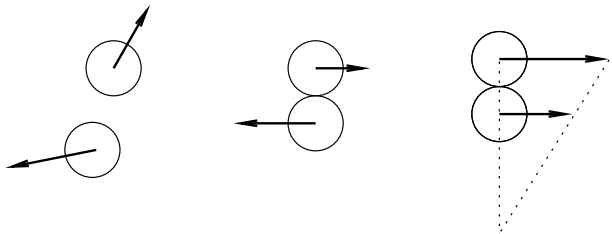}
\caption{Left: a broken contact. Center: a solid-like sheared contact. Right: a solid like stuck contact
(infinitesimal rigid rotation of a pair). The arrows indicate displacements.}
\end{center}
\end{figure}

The second type is a local shear motion.  In the local coordinates of one particle, it is either the motion of the second particle in the direction perpendicular to the line of centers, or an infinitesimal rotation (rolling). The corresponding contacts are called {\it sheared}. In our idealized model, friction in neglected, and any tangential force would lead to immediate separation of particles, because for disk-shaped particles, the contact surface is a point. We, however, still call these contacts solid-like, because in reality, these contacts are subject to friction forces, the contact surface has a positive area, and the particles in a sheared contact will not separate until the tangential force reaches the static friction threshold. In the case of rolling, the particles stay in contact and the pair is capable of bearing a compressive load.


Physically, vertices of the network can be realized as unit point masses and edges can be realized as elastic springs. Elastic force of the spring $(i,j)$ is determined by the pair potential $H(t_{ij})$, where $t_{ij}=({\bf u}^j-{\bf u}^i)\cdot{\bf q}^{ij}$. 
The potential is an important ingredient of our model, and therefore
we discuss it in detail. To motivate the choice
of $H$, we first 
recall the classical hard sphere potential
$H_{hs}$, which in our notation is defined by
\begin{equation}
\label{hs}
H_{hs}(t_{ij})=
\left\{
\begin{array}{cc}
\infty & \mbox{if}~t_{ij}<0,\\
0 & \mbox{if}~t_{ij}\geq 0.\\
\end{array}
\right.
\end{equation}
$H_{hs}$ models the following two options: 
(i) moving non-deformable (hard spheres) particles toward each other requires infinite energy
(a vertical line at $t_{ij}=0$), 
(ii) moving particles apart requires no energy.  Note that (\ref{hs}) already incorporates the constraints
(\ref{i3-intro}) by requiring infinite potential energy to violate the constraints.
\begin{figure}[htbp]
\label{potentials}
\begin{center}
\input{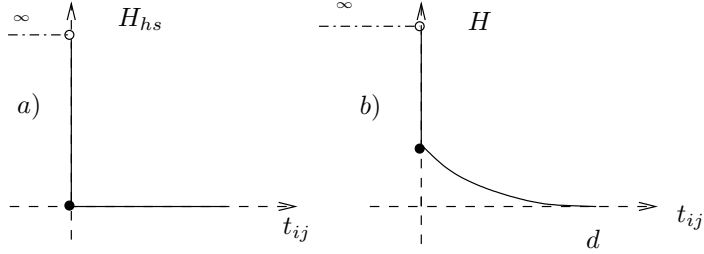}
\caption{
a) Hard sphere potential $H_{hs}$; b) The elastic potential $H(t_{ij})$ combines a vertical wall at $t_{ij}=0$ and a quadratic function with the vertex at $d$. }
\end{center}
\end{figure}

Elastic interaction between $D_i$ and $D_j$, together with constraints (\ref{i3-intro})
can be modeled by the following potential
\begin{equation}
\label{actual}
H(t_{ij}, d)=
\left\{
\begin{array}{cc}
\infty & \mbox{if}~t_{ij}<0,\\
\frac 12 Cd^{-3}(t_{ij}-d)^2 & \mbox{if}~t_{ij}\geq 0.\\
\end{array}
\right.
\end{equation}
\noindent
Here $d$ characterizes the cut-off distance of the potential, and $C$ determines
the magnitude of the pre-stress potential (the value of the potential when $t_{ij}=0$.
The potential (\ref{actual}) is shown on Fig. 2, together with the hard-sphere potential.
The formula (\ref{actual}) describes two options: (i) moving particles toward each other 
requires infinite energy; (ii) movement of particles apart from each other is caused
by finite, linear elastic force ${\bf f}^{ij}=-\partial H/\partial {\bf u}^i$; This force is repulsive for
small distances ($t_{ij}<d$), since $\partial H/\partial t^{ij}<0$. The magnitude of the force 
$f^{ij}(d)=\left|\frac{\partial H(t_{ij}, d)}{\partial t_{ij}}\right|=C\frac 12d^{-3}|t_{ij}-d|$. 
The magnitude of the force of pre-stress (or confining stress) is 
given by $\lim_{t_{ij}\to 0^+} f^{ij}=\frac 12 Cd^{-2}$, which tends to zero as $d\to\infty$
and $C$ is fixed.
Therfore, the effect of pre-stress is smaller for larger $d$. Further,
$H(t_{ij}, d)$ regularizes $H_{hs}$ in the following sense: if $d\geq d_0>0$, then
$\lim_{d\to\infty} H(t_{ij}, d)=0$ uniformly
on $(0, d_0]$. In the paper, we do not pass to this limit. Instead we choose $d$ sufficiently large and fix it, so
that $H$ is close to $H_{hs}$. Also, for technical simplicity, we set $C=1$ in (\ref{actual}), which corresponds
to an appropriate rescaling.

In reality, once the distance between $D_i$ and
$D_j$ is greater than the sum of their radii $a_i+a_j$, the pair interaction force
is zero. In our model, we still have a small repulsive force for all $a_i+a_j\leq t_{ij}\leq d$. So, the particles in our model
would continue accelerating away from each other even after separating. This favors separation of particles, and could lead to
increase in the number of broken contacts. Since our goal is estimating the number of solid-like contacts from below, this increase is acceptable.  In fact, our estimate holds 
for all $d$ sufficiently large.
We also mention that elastic contact force predicted by the classical Hertz theory
is a non-linear function of $t_{ij}$. Our model is chosen for simplicity, and can be viewed
as an approximation of Hertz theory, valid for sufficiently small displacements.

In general, the cut-off parameter $d$ may be different for different pairs of particles
in contact. Therefore, we define a pair interaction energy
\begin{equation}
\label{energy}
h(t_{ij}, d_{ij})= \frac 12d^{-3}(t_{ij}-d_{ij})^2,
\end{equation}
and consider 
\begin{equation}
\label{actual-gen}
H(t_{ij}, d_{ij})=
\left\{
\begin{array}{cc}
\infty & \mbox{if}~t_{ij}<0,\\
h(t_{ij}, d_{ij}) & \mbox{if}~t_{ij}\geq 0,\\
\end{array}
\right.
\end{equation}
where
\begin{equation}
\label{reg}
d_{ij}=d\delta_{ij},~~~~\delta_{ij}\in [1/2, 1].
\end{equation}
The formula (\ref{actual-gen}) is more general than  (\ref{actual}). 
The choice of $d_{ij}$ in (\ref{reg}) ensures that for all pairs $(i,j)$, the points where $h(t_{ij}, d_{ij})=0$
and located in the interval $[1/2d, d]$. The number 1/2 is of no particular significance. 
Any number in the interval $(0, 1)$ would work just as well. We only need all $d_{ij}$ to have the same sign and comparable magnitudes controlled by $d$.
Finally, the total elastic interaction energy $Q$ of the network is obtained by summing up $h(t_{ij}, d_{ij})$
over all pairs $(i,j)$ corresponding to the edges of the network.

The equilibrium state of a granular material corresponds to a minimum of $Q$ subject to the constraints
(\ref{i3-intro}) and the appropriate boundary conditions. Since the functional $Q$ is quadratic, and the
constraints and boundary conditions are linear, this is a quadratic programming problem, studied in optimization theory (e.g., \cite{Hestenes}). In the language of optimization theory, the solid-like contacts (\ref{solid-intro}) correspond to the so-called {\it active} constraints, while
the constraints corresponding to the broken contacts (\ref{broken}) are called {\it inactive}. The question
addressed in this paper concerns the total number and spatial distribution of each type of constraints in the energy-minimizing configuration of the network. It appears that no general results of this type are currently available in  optimization theory. The present study makes use of the geometric features of the contact graph, in particular its rigidity properties (see above), to investigate the energy minimizer.  We also note the connection between our constrained variational problem  and continuum variational inequalities
\cite{Duvaut}, \cite{Kikuchi}, \cite{Kinderlehrer}, \cite{Panagi}. Our problem can be viewed as a discrete
variational inequality.

The main result of the paper is Theorem \ref{main-theorem} in Section \ref{sect:optimality}.
There we consider a packing whose
contact graph in the reference configuration is a triangulation of a convex polygonal domain. The packing is deformed by imposing displacement boundary conditions at the packing boundary. The boundary conditions model the motion of rigid walls in engineering experiments. We prove that for generic contact graphs (precisely defined in Definition \ref{reg-treg}) and generic pre-stresses (corresponding to the choice of $\delta_{ij}$ in (\ref{energy}), 
the constrained energy minimizer for sufficiently large $d$ provides a packing with at least two solid-like contacts per each particle.  There also some (non-generic) choices of $\delta_{ij}$ for which theorem does not provide a definite conclusion.

The network of solid-like contacts is the load-bearing structure. The network of broken contacts can be associated with the micro-bands \cite{Kuhn1}, \cite{Kuhn2} that appear during small shear deformations. The result implies that no particle can lose contact with {\it all} of its neighbors, which eliminates ``micro-avalanches". Put another way, loss of structural integrity in dense packings is evolutionary rather than catastrophic, so that shearing with a small displacement will first lead to dilatation, during which the packing becomes more loose everywhere, and only then local avalanches may occur.

Another useful consequence of theorem \ref{main-theorem} is as follows. It provides a lower bound on the order parameter, recently introduced in \cite{AT, VTA} as one of the main ingredients of the new phenomenological theory of dense granular flows proposed by Aronson and Tsimring. The order parameter 
$\rho$ characterizes the phase transition from solid to fluidized state. To define $\rho$ at an arbitrary point ${\bf x}$ of $\Omega$, one begins by fixing a mesoscopic averaging volume $V$ of characteristic size $h$
(e.g., a disk of radius $h$ centered at  ${\bf x}$).  Then, all solid-like contacts within $V$ should be counted. Next, the obtained number $n_s$ of solid-like contacts is divided by the number $n$ of all contacts within $V$
to obtain $\rho$. So, $\rho$ is a mesoscopic average, which in general depends on $h$.
In many systems, such as periodic elastic composites, the results of mesoscopic averaging is practically independent of $h$ for $h$ larger than a certain characteristic length. For disordered granular materials this
is not necessarily true. Therefore, a rigorous mathematical theory may require study of a family of order parameters parametrized by $h$. In Section 6, we define such a family of order
parameters using the notion of a $k$-neighborhood of a vertex of
$\Gamma$ as a discrete analogue of $V$.  Specifying an integer $k$ in our definition corresponds to
choosing $h$ in the continuous case. 


The paper is organized as follows. In Sect. 2 we formulate the main constrained minimization problem.
The problem contains two types of constraints: impenetrability constraints (\ref{i3-intro}) and the
boundary constraints (see (\ref{b1}), (\ref{b1-bis})), corresponding to the external boundary conditions.  Elimination of these boundary constraints
leads to a reduced minimization problem. In Sect. 3 we recall some facts concerning first-order rigidity of graphs. In Sect. 4 we show existence of a unique minimizer of the reduced
problem. Optimality conditions for the reduced problem are stated and analyzed in Sect. 5, where we also state and prove the main theorem \ref{main-theorem}. In Sect. 6 we introduce a definition of the order parameter in the spirit of \cite{VTA} and give a lower bound on the
order parameter that follows from the main theorem. Finally, conclusions are provided in Sect.
7.



\section{Formulation of the problem}
\subsection{Elastic interactions with impenetrability constraints}
\label{sect-constraints}
In 2D, consider a packing of spheres $D_i$ of radii
$a_i$ with centers ${\bf x}^i, i=1,2,\ldots, N$. (All vectors in this
paper are column vectors and we use superscript `{\sf T}' to indicate
transposition). The packing fills a bounded region. After an infinitesimal motion, the position of the center of $D_i$ is 
${\bf y}^i$. We write
${\bf y}^i={\bf x}^i+{\bf u}^i$ where ${\bf u}^i$ are displacements.
The vertices ${\bf x}^i, {\bf x}^j$ are connected by an edge if and only if $D_i, D_j$ are in contact. In this case, we call ${\bf x}^i$ and ${\bf x}^j$ neighbors. We denote 
by 
$
{\mathcal N}_i
$
the set of $j\in \{1, 2,\ldots, N\}$ such that ${\bf x}^j$ is a neighbor of ${\bf x}^i$.
Orientation of contacts (equivalently, edges) is prescribed by the unit vectors
\begin{equation}
\label{orient}
{\bf q}^{ij}=\frac{{\bf x}^j-{\bf x}^i}{|{\bf x}^j-{\bf x}^i|}.
\end{equation}
The vertices ${\bf x}^i$ and edges $(i,j)$ define the contact graph
$\Gamma$.  
Let $E$ denote the number of edges of $\Gamma$.
The edge set ${\mathcal E}$ of $\Gamma$ is given by
$\{(i,j): j\in{\mathcal N}_i, i =1,2, \ldots, N\}$. To each edge $(i,j)$
we can associate a pair potential energy $h(t_{ij}, d_{ij})$ defined in (\ref{energy}).
Summing up all these energies we obtain the total 
elastic interaction energy of the network. It is a quadratic form
\begin{equation}
\label{q-basic}
Q({\bf u}^1,{\bf u}^2, \ldots {\bf u}^N)=\sum_{i=1}^N
\sum_{j\in {\mathcal N}_i} h(t_{ij}, d_{ij})=
\frac 12 d^{-3}\sum_{i=1}^N
\sum_{j\in {\mathcal N}_i}
\left (({\bf u}^j-{\bf u}^i)\cdot{\bf q}^{ij}-d_{ij}\right)^2,
\end{equation}
on the displacements ${\bf u}^i, i=1,2,\ldots, N$.
In (\ref{q-basic}), $d, d_{ij}$ are parameters specified by (\ref{actual-gen}), (\ref{reg}). 

Our objective is to determine the displacements ${\bf u^i}$, $i=1,2, \ldots,
N$ so that the the energy functional $Q$ is minimized subject to two
types of constraints. The first type of constraints consists of
{\em linearized impenetrability constraints}. These are obtained by formally
linearizing the condition
that the distance between two spheres in contact cannot decrease. Consider
two spheres $D_i, D_j$ in contact. In the reference configuration,
\begin{equation}
\label{i1}
|{\bf x}^i-{\bf x}^j|=a_i+a_j.
\end{equation}
Assuming that $D_i, D_j$ cannot overlap, we have
\begin{equation}
\label{i2}
|{\bf y}^i-{\bf y}^j|\geq a_i+a_j.
\end{equation}
These are the impenetrability constraints.
We linearize (\ref{i2}) by writing
\begin{eqnarray*}
|{\bf y}^i-{\bf y}^j|^2 & = & |{\bf x}^i-{\bf x}^j+{\bf u}^i-{\bf u}^j|^2\\
          & = &|{\bf x}^i -{\bf x}^j|^2 +
            2({\bf x}^i - {\bf x}^j) \cdot ({\bf u}^i - {\bf u}^j) +
            |{\bf u}^i - {\bf u}^j|^2 \\
          & = & (a_i + a_j)^2 + 
            2({\bf x}^i - {\bf x}^j) \cdot ({\bf u}^i - {\bf u}^j) +
            |{\bf u}^i - {\bf u}^j|^2. 
\end{eqnarray*}
Now for for ``small'' $|{\bf u}^i - {\bf u}^j|$ we can neglect quadratic term
$|{\bf u}^i - {\bf u}^j|^2$, and (\ref{i2}) yields $2({\bf x}^i - {\bf x}^j)\cdot ({\bf u}^i -{\bf u}^j) \geq
0$, which in turn is equivalent to $({\bf u}^j - {\bf u}^i)\cdot
{\bf q}^{ij} \geq 0$ where $ {\bf q}^{ij}$ is as defined in
(\ref{orient}).
Therefore, the first set of constraints we impose on the displacements
${\bf u}^1, {\bf u}^2, \ldots, {\bf u}^N$ is
\begin{equation}
\label{i3}
({\bf u}^j-{\bf u}^i)\cdot {\bf q}^{ij}\geq 0, ~~~j\in {\mathcal N}_i,
~~i=1,2,\ldots, N.
\end{equation}
The second type of constraints corresponds to the boundary conditions.
Particles located at the packing boundary have {\em prescribed} displacements.
In the sequel we refer to these particles as {\it boundary particles}. The corresponding vertices of $\Gamma$ are called {\it boundary vertices}.
Other particles are referred to as  {\it interior}, or sometimes, {\it free}, 
and the corresponding vertices of $\Gamma$ as {\it interior vertices}. 

All boundary particles are divided into several groups,
numbered $1,2,\ldots,M$. Let $I_m$ denote the set of indices
of the particles in group $m$ for $m=1,2, \ldots, M$.
Each sphere in a certain group is in contact with at least one
other sphere from the same group. 
Each group moves as a single rigid body. We assume that the prescribed boundary displacements are of the form
\begin{equation}
\label{b1}
{\bf u}^i={\bf R}^m({\bf x}^i),~~~i\in I_m, ~~m=1,2,\ldots, M,
\end{equation}
where
\begin{equation}
\label{b1-bis}
{\bf R}^m({\bf x}^i)={\bf c}^m+\alpha^m K({\bf x}^i-{\bf x}^{\star,m}),
~~~i\in I_m, ~~m=1,2, \ldots, M,
\end{equation}
and ${\bf c}^m, {\bf x}^{\star,m}$ are given vectors, $\alpha^m$ is a
given scalar, and
$K$ is the matrix denoting clockwise rotation by $\pi/2$.
The functions ${\bf R}^m$ are called {\em infinitesimal rigid
displacements},
parametrized by a scalar $\alpha^m$, and
vectors ${\bf c}^m$ and ${\bf x}^{\star,m}$. We refer the
reader to Sect. \ref{sect:rigidity} for more details
on rigid displacements.

Our description above leads to the\\

\vspace*{0.5\baselineskip}
\noindent
\underline{\bf Main problem:} 
\begin{eqnarray}
\label{mpr1}
\mbox{minimize}&~~&  Q({\bf u}^1, {\bf u}^2, \ldots, {\bf u}^N)\\
\label{mpr2}
\mbox{subject to}&~~& \mbox{linearized impenetrability
                            constraints (\ref{i3})}\\
\label{mpr3}
                 &~~& \mbox{and boundary conditions (\ref{b1}).}
\end{eqnarray}

\subsection{Feasible region}
Let us define the configuration space $U$. Points of this space
are denoted by ${\bf U}=(({\bf u}^1)\trp, ({\bf u}^2)\trp, \ldots,
({\bf u}^N)\trp)\trp$.

\noindent
{\bf Remark}. To avoid this heavy notation, we simply write
$$
{\bf U}=({\bf u}^1, {\bf u}^2, \ldots, {\bf u}^N),
$$
when no confusion can occur.

Dimension of $U$ is $2N$.  
{\it Feasible region} ${\mathcal F}$ is the subset of $U$
in which all the constraints (\ref{i3}) and (\ref{b1})
are satisfied. The points satisfying (\ref{i3}) form a polyhedral
(not necessarily bounded) region. The boundary of this region
consists of parts of the hyperplanes (subspaces of dimension $2N-1$) defined
by 
\begin{equation}
\label{i3-eq}
({\bf u}^j-{\bf u}^i)\cdot {\bf q}^{ij}=0,~~~~j\in{\mathcal N}_i,~i=1, 2, \ldots, N.
\end{equation}
Because of the close relation to rigidity, we refer to (\ref{i3-eq}) as
$R$-{\it equations}. 
Equations
(\ref{b1}) define
$M$ planes $S_m, m=1, \ldots, M$. Dimensions
of $S_m$ depend on the number of the boundary particles in the $m$-th group. 

For each point of ${\bf U}\in {\mathcal F}$, 
some of the constraints (\ref{i3}) are satisfied as equations. These constraints
are called {\it active}. The corresponding edges of the contact graph $\Gamma$ 
are called active as well. 
The rest of (\ref{i3}) are satisfied as strict inequalities. These are {\it inactive}
constraints (respectively, edges).

\subsection{Elimination of constraints corresponding to boundary conditions}
\label{sect:elimination}
The quadratic form $Q$ in (\ref{q-basic}) can be written in a convenient
form in terms of a certain matrix $R^r$. To define $R^r$, we
index the edges of $\Gamma$ by $l$,
$l=1,2, \ldots, E$. Let $(i_l, j_l) \in {\mathcal E}$ be the edge
of $\Gamma$ corresponding to $l$ for $l=1,2,\ldots, E$. 
Let $R^r$ be the $E\times 2N$ matrix whose $l$-th row is defined by
\begin{equation}
\label{row}
R^r_{lm} = \left\lbrace\begin{array}{ll}
             +({\bf q}^{i_lj_l})_1 &\mbox{if $m=2(j_l-1)+1$}\\
             +({\bf q}^{i_lj_l})_2 &\mbox{if $m=2(j_l-1)+2$}\\
             -({\bf q}^{i_lj_l})_1 &\mbox{if $m=2(i_l-1)+1$}\\
             -({\bf q}^{i_lj_l})_2 &\mbox{if $m=2(i_l-1)+2$}\\
             0                     &\mbox{otherwise}
             \end{array}\right .
\end{equation}
for $l=1,2, \ldots, E$.

\noindent
{\bf Remarks.} 
1. $R^r$ is the (first-order) rigidity matrix,
a well known object in geometric rigidity theory (see e.g. \cite{CoW, W}).

\noindent
2. Consider
vertices ${\bf x}^{i_l}, {\bf x}^{j_l}$ and the edge $l$ connecting them.
The corresponding row ${\bf r}^l$ of $R^r$ has $2N$ entries. We can view ${\bf r}^l$
as a string of $N$ pairs of numbers, the first pair corresponding to ${\bf x}^1$,
the second to ${\bf x}^2$ and so on. For simplicity, we shall call a pair
of entries corresponding to a particular vertex ${\bf x}^i$
a {\it place corresponding to ${\bf x}^i$}.

Then we can interpret equation (\ref{row}) as follows. A row ${\bf r}^l$ has zeros at all places, except two. The non-zero entries are $-{\bf q}^{i_l, j_l}$, written as a two-dimensional row vector
at the place corresponding to ${\bf x}^{i_l}$; and ${\bf q}^{i_l, j_l}$,
written as a two-dimensional row
at the place corresponding to ${\bf x}^{j_l}$.

\noindent
3. A row of $R^r$ corresponds to an edge of $\Gamma$. Therefore it is natural to call a row active (respectively, inactive)
if a corresponding edge is active (respectively, inactive).

Now define the vector ${\bf d}\in\reals^E$ by
\begin{equation}
\label{d}
{\bf d}=-(d_{i_1j_1},~d_{i_2j_2},~ \ldots, ~d_{i_Ej_E}), 
\end{equation}
where $d_{i_l, j_l}$ are chosen according to (\ref{reg}).               
With these notations the quadratic form $Q$ in (\ref{q-basic})
can be written as
\begin{equation}
\label{q}
Q({\bf U})=d^{-3}\frac 12 (R^r{\bf U}+{\bf d})\cdot( R^r{\bf U}+{\bf d}). 
\end{equation}

We now eliminate the boundary conditions (\ref{b1}) from the main
problem (\ref{mpr1}, \ref{mpr2}, \ref{mpr3}).
Let $N_b= \sum_{m=1}^M \mbox{card}(I_m)$.
Then the equations (\ref{b1}) simply state that the $2N_b$
components of ${\bf U}$ corresponding to the $N_b$
boundary vertices have prescribed
displacements. Without loss of generality assume
that the {\em last} $2N_b$ components of $U$ correspond
to the boundary vertices.
Let us partition ${\bf U}$ as
\begin{equation}
\label{split1} 
{\bf U}=[~{\bf z}~|~{\bf w}~],
\end{equation} 
where ${\bf z}=({\bf u}^1, {\bf u}^2, \ldots,
{\bf u}^{N-N_b})$ corresponds
to displacement vectors of interior vertices, and 
${\bf w}= ({\bf u}^{N-N_b+1}, {\bf u}^{N-N_b+2}, \ldots,
{\bf u}^{N})$ corresponds to the displacements of the
boundary vertices.
The equality constraint (\ref{b1}) is now simply
\begin{equation}
\label{b3}
{\bf w}={\bf g}
\end{equation} 
where $g\in \reals^{2N_b}$ is the vector of displacements prescribed by
the right-hand-sides of (\ref{b1}).
The matrix $R^r$ can be partitioned similarly to (\ref{split1}):
\begin{equation}
\label{split2}
R^r=[~R~|~R^b~],
\end{equation}
where dimensions of $R$ and $R^b$ are $E\times 2(N-N_b)$ and $E\times
2N_b$, respectively. 
Denote 
\begin{equation}
\label{a}
{\bf a}=R^b{\bf g}.
\end{equation}
Using (\ref{split1})--(\ref{a}) in (\ref{q}) and in
(\ref{i3}) we can reduce the main problem
(\ref{mpr1}, \ref{mpr2}, \ref{mpr3})
to\\

\vspace*{0.5\baselineskip}
\noindent
\underline{{\bf Reduced problem:}}
\begin{eqnarray}
\label{pr1}
\mbox{minimize}&~~&F({\bf z})=\frac 12 d^{-3}\left(R{\bf z}+{\bf a}+{\bf d}\right)
\cdot\left( R{\bf z}+{\bf a}+{\bf d}\right)\\
\label{pr2}
\mbox{subject to}&~~&R{\bf z}+{\bf a}\geq 0.
\end{eqnarray}
\vspace*{0.5\baselineskip}
The minimization in (\ref{pr1}) is taken over all ${\bf z}\in
\reals^{(N-N_b)}$.

\section{First-order rigidity}
\label{sect:rigidity}
A rigid motion is a composition of a translation and rotation:
\begin{equation}
\label{r1}
{\bf y}({\bf x})={\bf c}+{\bf x}^\star+O({\bf x}-{\bf x}^\star),
\end{equation}
where $O$ is an orthogonal (rotation) matrix, ${\bf c}$ is a translation vector,
${\bf x}^\star$ is a center of rotation. If $O$ is close to identity $I$ (infinitesimally small rotation), then
$$
O\approx I+A,
$$
where $A$ is a skew matrix ($a_{ij}=-a_{ji}$).

Suppose that in a two-dimensional rigid motion, the rotation angle $\alpha$ is close to zero.  
Then
$$
O=
\left(
\begin{array}{cc}
\cos \alpha & \sin \alpha \\
-\sin \alpha & \cos \alpha\\
\end{array}
\right)\approx
\left(
\begin{array}{cc}
1 & 0 \\
0 & 1\\
\end{array}
\right)
+
\alpha
\left(
\begin{array}{cc}
0 & 1 \\
-1 & 0\\
\end{array}
\right)=I+\alpha K,
$$
where 
$$
K=\left(
\begin{array}{cc}
0 & 1 \\
-1 & 0\\
\end{array}
\right)
$$ 
is a clockwise rotation by $\pi/2$. In that case, (\ref{r1}) becomes
\begin{equation}
\label{r2}
{\bf y}({\bf x})={\bf c}+{\bf x}+\alpha K({\bf x}-{\bf x}^\star)=
\left(
\begin{array}{c}
c_1+x^\star_1-\alpha (x-x^\star)_2\\
c_2+x^\star_2+\alpha (x-x^\star)_1\\
\end {array}
\right).
\end{equation}

Let ${\bf u}={\bf y}({\bf x})-{\bf x}$ denote the displacement. We can write
(\ref{r2}) as
\begin{equation}
\label{r3}
{\bf u}({\bf x})={\bf c}+\alpha K({\bf x}-{\bf x}^\star).
\end{equation}

\begin{definition}
We call (\ref{r3}) an infinitesimal rigid displacements in 2D.
\end{definition}

Next, let $G$ be a graph. 
Consider all motions of vertices of $G$ that preserve the lengths of the edges.
If the only such motions are the rigid body motions of the whole graph, 
then the graph is called {\it rigid}.
A graph $\Gamma$ is {\it first-order rigid} \cite{W} if all solutions
of the $R$-system (\ref{i3-eq})
are infinitesimally rigid displacements.  
In addition, $\Gamma$ is {\it independent}
if
the rows of the rigidity matrix $R^r$ are linearly independent. Graphs that are both
first-order rigid and independent are called {\it isostatic} (\cite{W}). Intuitively, an isostatic graph is minimally rigid, that is removing any edge results in loss of rigidity. Another notion of rigidity
is {\it generic rigidity}, see \cite{W}. According to Thm. 49.1.7 from \cite{W},
generic rigidity for a neighborhood in a configuration space is equivalent to the first-order rigidity for some specific configuration in that neighborhood.  

With the definition of $R^r$ and ${\bf U}$ in Sect. \ref{sect:elimination}
the system (\ref{i3-eq}) can be written as
\begin{equation}
\label{r-sys-a}
R^r{\bf U}=0.
\end{equation}

Note the connection between $R$-system and constraints, 
as well as the functional of the main problem
(\ref{mpr1},{\ref{mpr2},{\ref{mpr3}).

The following definition (see \cite{W} for a $d$-dimensional definition) is useful for verifying rigidity of graphs.
\begin{definition}
\label{henneberg}
For a graph $\Gamma$, the Henneberg 2-construction in 2D is a sequence of graphs $G_1, G_2, \ldots, G_n$ such that:\newline
\noindent
(i) $G_{k+1}$ is obtained from $G_k$ by either vertex addition (attaching a new vertex by 2 edges); or edge splitting (replacing and edge from $G_k$
with a new vertex joined to its ends and to 1 other vertex);\newline
\noindent
(ii) $G_k$ is a complete graph on $k$ vertices, and $G_n=\Gamma$.
\end{definition}
The following result
is stated in (\cite{W}, thm. 49.1.13):
\begin{theorem}
\label{hen-thm}
If a graph $G\subset {\reals}^2$ is obtained by a Henneberg $2$-construction, then $G$ is 
generically isostatic. 
\end{theorem}
\noindent
{\bf Remark 1}.
In veiw of the definitions above, Theorem \ref{hen-thm} implies that a graph obtained by Henneberg $2$-construction is first-order rigid and independent (minimally rigid).

In the present case, the rows of the rigidity matrix $R^r$ are not linearly independent, but the row rank is maximal. This means that we typically have more edges than needed to ensure rigidity of $\Gamma$. In this situation, the following theorem (thm. 49.1.14 from \cite{W}) is useful.
\begin{theorem}
\label{gluing}
If two graphs $G_1$ and $G_2$ are generically rigid planar graphs sharing at least $2$ vertices,
then the graph $G$ obtained by combining all vertices and edges of $G_1, G_2$ is generically rigid.
\end{theorem}

\noindent
{\bf Remark 2}.  Because of the relation between generic rigidity and first-order rigidity, Theorem
\ref{gluing} implies that combining first-order rigid graphs $G_1$, $G_2$ as in this Theorem
yields a first-order rigid graph $G$. We shall use Theorem \ref{gluing} to obtain first-order rigidity of
triangulations. Indeed, one triangle $G_1$ is first-order rigid. Adding another triangle $G_2$ so that
$G_1$ and $G_2$ share an edge, yields a first-order rigid graph. Then we can proceed sequentially.
Given  a first-order rigid triangulation $G_k$, we construct $G_{k+1}$ by combining $G_k$ with a triangle. This new triangle either shares two vertices with $G_k$, or all three vertices. In the first case, $G_{k+1}$ would have
one more vertex and two more edges than $G_k$. In the second case, $G_{k+1}$ has the same number of vertices as $G_k$, and one more edge. By theorem \ref{gluing}, any planar triangulation obtained by this sequential procedure is first-order rigid.


\section{Existence and uniqueness of minimizers of the reduced problem}
Let $\Omega$ be a bounded connected domain in $\reals^2$ with a polygonal boundary.
First we show that, under certain assumptions on geometry of $\Gamma$, the matrix $R$ has full column rank.

We shall say that $\Gamma$ is a triangulation if edges of $\Gamma$ partition $\Omega$ into a disjoint
union of triangles.  

Let $X$ be a set of interior vertices of $\Gamma$, containing at least two elements. Consider a graph
$G_X\subset \Gamma$ defined as follows. Vertices of $G_X$ are all elements of $X$. Edges of 
$G_X$ are those edges of $\Gamma$ that join two vertices from $X$. We also assume that $X$ is chosen
so that $G_X$ is a connected graph.
\begin{definition}
\label{solid}
The contact graph $\Gamma$ is {\it cell-connected} if for each $G_X\subset \Gamma$ as above, there
exist two vertices ${\bf x}^1, {\bf x}^2$ in $G_X$, and two interior vertices $\hat{{\bf x}}^1, \hat{{\bf x}}^2$ in $\Gamma\backslash G_X$, such that the quadrilateral with vertices 
${\bf x}^1, {\bf x}^2, \hat{{\bf x}}^1, \hat{{\bf x}}^2$ is a union of two {\it adjacent} triangles of $\Gamma$.
\end{definition}
An example illustrating the definition if shown in Fig. 3.
\begin{figure}[]
\label{non-connected}
\begin{center}
\input{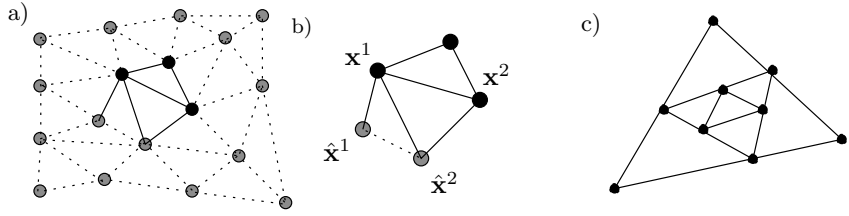}
\caption{a) --A cell-connected graph. Solid dots indicate vertices of $G_X$;
b) -- The cell connection with three edges and four vertices is shown separately; 
c)-- A triangulation that is not cell connected. Here the vertices of the small triangle inside
are connected to other vertices by pairs of collinear edges.}
\end{center}
\end{figure}

\begin{definition}
\label{reg-treg}
We call $\Gamma$ a {\it regular triangulation} if \newline
(i) $\Gamma$ can be obtained by sequential addition of triangles in the way described in the Remark 2;\newline
\noindent
(ii) every interior vertex, connected to a boundary vertex,  is also connected to at least one other
boundary vertex, and the corresponding edges are non-collinear;\newline
\noindent
(iii) $\Gamma$ is cell-connected.
\end{definition}

\noindent
{\bf Remark. } 
Note that i) in Definition \ref{reg-treg} implies that $\Gamma$ is first-order rigid.
The property ii) states that every edge connecting a boundary vertex 
with an interior vertex must be a part of the boundary of a triangle, containing two boundary vertices and one  interior vertex. 

Informally, Definition  \ref{solid} (or (iii) in Definition \ref{reg-treg} )
can be interpreted as a strong connectivity property. According to Definition \ref{solid}, a generic connected subgraph
$G_X$ is connected to the ``rest of $\Gamma$" not just by an edge, but by a "more robust" cell structure that consists of three edges, with one edge bracing the other two, (see Fig 3.). Note also that either ${\bf x}^1$ or ${\bf x}^2$ are connected to the vertices of $\Gamma\backslash G_X$ by a pair of non-collinear edges. This observation is important in the proof of Proposition \ref{full-rank} below. It is not difficult to see
that a periodic triangular planar graph satisfies iii). However, there are triangulations with a mean coordination number four that do not satisfy iii). An example of such a graph is shown in Fig. 3 c).


\begin{proposition}
\label{full-rank}
Suppose that $\Gamma$ is a regular triangulation. Then ${\mbox rank}~R=2(N-N_b)$.
\end{proposition}

\noindent
{\it Proof}. Consider a subgraph $\Gamma_{max}\subset \Gamma$ constructed inductively as follows. 
Begin with $\Gamma_1$ that consists of all boundary vertices. On the next step, add an interior vertices connected to $\Gamma_1$ by two or more non-collinear edges. Also, add exactly two non-collinear edges that connect this vertex to $\Gamma_1$. Call
the resulting graph $\Gamma_2$. Generally, given $\Gamma_k$, $k\geq 2$, define
$\Gamma_{k+1}=\Gamma_k\cup S_k$, where $S_k$ consists of an interior vertex ${\bf x}^k$, not contained in
$\Gamma_k$ but connected to $\Gamma_k$ by at least two non-collinear edges, together with
a pair of non-collinear edges connecting ${\bf x}^k$ to $\Gamma_k$. Since the graph $\Gamma$ has a finite number of vertices, the
process terminates after a finite number of steps. The resulting graph is $\Gamma_{max}$. The construction
is illustrated in  Fig. 
\begin{figure}[]
\label{non-connected}
\begin{center}
\input{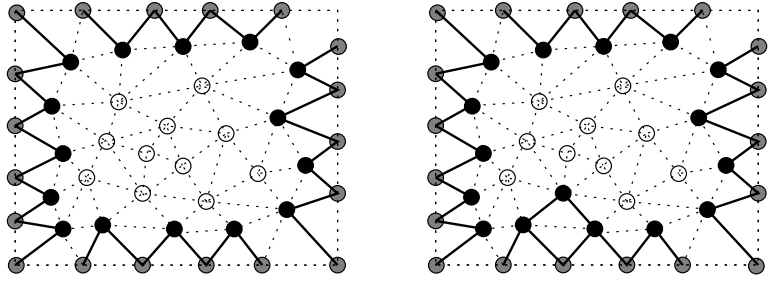}
\caption{The subgraphs $\Gamma_1$ (left) and $\Gamma_2$ (right). The boundary vertices are shown
in gray. The interior vertices of $\Gamma_1$ and $\Gamma_2$ are shown in black. The edges of
$\Gamma_1$, $\Gamma_2$ are shown by solid lines. The other vertices of $\Gamma$ are represented by
the unfilled circles. The edges of $\Gamma$ not included into $\Gamma_1, \Gamma_2$ are represented by the dotted lines.}
\end{center}
\end{figure}
\begin{figure}[]
\label{non-connected}
\begin{center}
\input{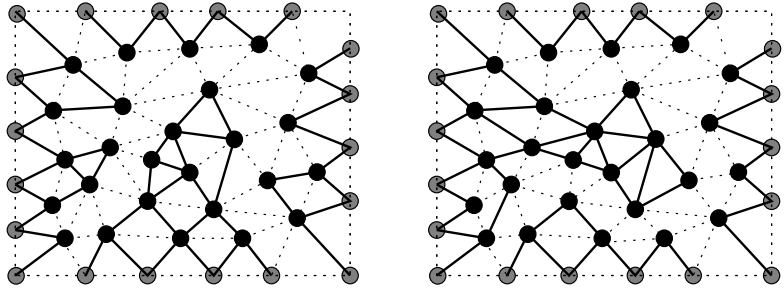}
\caption{Two different subgraphs $\Gamma_{max}$ constructed inductively. Both subgraphs
are constructed starting with $\Gamma_2$ in Fig. 4.}
\end{center}
\end{figure}
We claim that $\Gamma_{max}$ contains all vertices of $\Gamma$. To obtain
a contradiction, suppose that there are vertices not included into $\Gamma_{max}$.
Denote the set of these vertices by $X$, and denote by $G_X$ the graph formed by vertices in
$X$ and all edges of $\Gamma$ that connect these vertices. Let $G_X^c$ be any connected component
of $G_X$. If $G_X^c$ is a single point ${\bf x}_g$, then, since $\Gamma$ is a triangulation, there must be at least three edges incident at ${\bf x}_g$, and at least two of these edges must be non-collinear. 
Thus  ${\bf x}_g$ must be in $\Gamma_{max}$, which gives a contradiction.
Next suppose that $G_X^c$ contains two or more vertices. By Definition \ref{solid}, (see also the Remark following that Definition), there must be a pair of vertices of ${\bf x}^1, {\bf x}^2$ in $G_X^c$ connected to two vertices in $\Gamma\backslash G_X$ by three edges, and at least one of ${\bf x}^1, {\bf x}^2$ must
be connected to $\Gamma\backslash G_X$ by two non-collinear edges. Denote this vertex by ${\bf x}_g$. It must be included into $\Gamma_{max}$ which gives a contradiction and proves the claim.

Next, we claim that the number of edges in $\Gamma_{max}$ is $2(N-N_b)$. Indeed, each interior
vertex in $\Gamma_2$ has exactly two non-collinear edges incident at it. Then, on each of the next steps, we add an interior vertex together with two non-collinear edges incident at it. Since
the number of free vertices in $\Gamma_{max}$ is $N-N_b$, the claim is proved.

Finally, we claim that the rows of $R$ corresponding to the edges of $\Gamma_{max}$ are linearly independent. Let the matrix of these rows be denoted by
$R_{max}$. This is a square $2(N-N_b)$ matrix. We claim that an appropriate
row-reduction reduces $R_{max}$ to a matrix $R^\prime_{max}$ that has
block-diagonal form: for each
vertex ${\bf x}^i$ there are exactly two rows ${\bf r}^{i,1}, {\bf r}^{i,2}$ of $R^\prime_{max}$,
and two linearly independent unit vectors ${\bf q}^{ij_1}$ and ${\bf q}^{ij_2}$,
such that ${\bf r}^{i,1}$ (${\bf r}^{i,2}$) contains
${\bf q}^{ij_1}$ (${\bf q}^{ij_2}$) at a place
corresponding to ${\bf x}^i$ while all other entries in these rows are zero.

To see this, consider first a ``basic unit" of $\Gamma_2$: an interior vertex ${\bf x}^i$ and two non-collinear edges
incident at it. Let the corresponding unit vectors be ${\bf q}^{i1}$, ${\bf q}^{i2}$. Recall that these edges connect ${\bf x}^i$ to two boundary vertices. Consequently,
the rows ${\bf r}^1, {\bf r}^2$ corresponding to the above pair of edges have zeros
at all places, except two places corresponding to ${\bf x}^i$. 
The non-zero entries of ${\bf r}^1$ (${\bf r}^2$) are two components of ${\bf q}^{i1}$
(${\bf q}^{i2}$).
Since ${\bf q}^{i1}$, ${\bf q}^{i2}$ are linearly independent, so are
${\bf r}^1, {\bf r}^2$.  Furthermore, linear combinations
of ${\bf r}^1, {\bf r}^2$ can be used to eliminate non-zero entries in other rows.
By adding an appropriate linear combination of ${\bf r}^1, {\bf r}^2$ to a row with
some unit vector ${\bf q}^{ij}$ at a place corresponding to ${\bf x}^i$, we can obtain
zeros at this place. Hence, by using ${\bf r}^1, {\bf r}^2$ as pivots in Gaussian elimination
we can obtain rows whose only non-zero entries are at a place corresponding
to a vertex that was added to $\Gamma_2$ on the next step of the iteration. These rows, in turn, can be used as pivots. Continuing
with row reduction, we can eventually reduce all rows of $R_{max}$ to this form
(this follows from the fact that $\Gamma_{max}$ contains all vertices of $\Gamma$). 
The proposition is proved.

\noindent
{\bf Remark}. It is interesting to compare $R$ and the rigidity matrix $R^r$. It is 
well-known that for a first-order rigid graph,  the null space of $R^r$ is non-trivial
and consists of infinitesimal rigid displacements. The proposition above shows that
the null space of $R$ is trivial. The main difference in structure between these two 
matrices is that $R$ contains special rows, that might be called {\it broken}. These
rows correspond to edges connecting an interior vertex to a boundary one. 
A typical row of $R^r$ has four non-zero entries, while each broken row
has only two.  These entries occur at a place corresponding to an interior vertex.
If an interior vertex is connected to two boundary vertices, then the regular
triangulation property of $\Gamma$ ensures that the corresponding broken rows are non-collinear, and can
be used in the sequential Gaussian elimination, as done in the proof.

\begin{proposition}
\label{existence}
Consider the problem (\ref{pr1},\ref{pr2}). Suppose that $\Gamma$ is a regular triangulation,
the feasible set of (\ref{pr1},\ref{pr2}) is non-empty,
and that the unconstrained minimizer of $F({\bf z})$ is not feasible. 
Then
the problem (\ref{pr1},\ref{pr2}) admits a unique minimizer that is a point
on the boundary of its feasible set.
\end{proposition}

\noindent
{\it Proof.} \vspace{0.1cm}
The problem (\ref{pr1},\ref{pr2}) has a feasible
point $\bar {\bf z}$. Then the problem (\ref{pr1},\ref{pr2}) has a unique
minimizer ${\bf z}^*$ as we now demonstrate.
Let the set ${\mathcal L}(\bar {\bf z})$
be defined by
\begin{equation*}
{\mathcal L}(\bar {\bf z}) = \{{\bf z}\in\reals^{2(N-N_b)}:
F({\bf z}) \leq F(\bar {\bf z})\}.
\end{equation*}
By Proposition \ref{full-rank} the matrix $R$ has full column rank. Therefore, the set
${\mathcal L}(\bar {\bf z})$ is an ellipsoid, which is a closed,
bounded, convex set.
The set
${\mathcal H}$ of ${\bf z} \in \reals^{2(N-N_b)}$ satisfying the constraint
(\ref{pr2})
is a closed half space. Therefore, ${\mathcal L}(\bar {\bf z}) \cap {\mathcal H}$
is a nonempty, closed, bounded convex set. Indeed, the reduced problem
(\ref{pr1},\ref{pr2}) is equivalent to the problem
\begin{eqnarray}
\label{preq1}
\mbox{minimize}&~~&F({\bf z})\\
\label{preq2}
\mbox{subject to}&~~&{\bf z} \in {\mathcal L}(\bar {\bf z})\cap{\mathcal H}.
\end{eqnarray}
Now by the continuity of $F$ and the compactness of
${\mathcal L}(\bar {\bf z})\cap{\mathcal H}$ we see that the problem
(\ref{preq1},{\ref{preq2}) and hence the reduced problem
(\ref{pr1},\ref{pr2}) has a minimizer ${\bf z}^*$. Since $R$ has full
column rank, $R\trp R$ is positive definite. The positive definiteness
of the matrix $R\trp R$ implies that $F$ is strictly convex
on ${\mathcal L}(\bar {\bf z})\cap {\mathcal H}$, from which we conclude that
${\bf z}^*$ must be unique.

Now if $R{\bf z}^\star+{\bf a}>0$, then ${\bf z}^\star$ must be the
unconstrained minimizer of $F$. This contradicts the assumption that
the unconstrained
minimizer is not feasible. Therefore, some components of
$R{\bf z}^\star+{\bf a}$ must be zero, and thus ${\bf z}^\star$ must be on the boundary
of the feasible region.

\section{Optimality conditions for the reduced problem}
\label{sect:optimality}
The reduced problem (\ref{pr1},\ref{pr2}) has a convex objective
function and linear constraints. For such problems 
it is possible to state optimality conditions 
that  are both necessary and sufficient \cite{NW}.
To be specific define the Lagrangian $L$ for the problem
({\ref{pr1},\ref{pr2}) by
\begin{equation}
\label{L}
L({\bf z}, \boldsymbol \lambda)=\frac 12 d^{-3}(R{\bf z}+{\bf a}+{\bf d})\cdot(R{\bf z}+{\bf a}+{\bf d})-
d^{-3}\boldsymbol \lambda\cdot (R{\bf z}+{\bf a}).
\end{equation}

Then ${\bf z}^*$ solves problem (\ref{pr1},\ref{pr2}) if and only if
there exists $\boldsymbol\lambda^*$ such that
${\bf z}={\bf z}^*$ and $\boldsymbol\lambda =\boldsymbol\lambda^*$
satisfy the Karush, Kuhn, Tucker (KKT) conditions
\begin{eqnarray*}
&& \nabla_{z} L ({\bf z},\boldsymbol\lambda) = 0\\
&& R{\bf z} +{\bf a}\geq 0\\
&& \boldsymbol\lambda\cdot (R{\bf z}+{\bf a})=0\\
&& \boldsymbol\lambda\geq 0
\end{eqnarray*}
or equivalently
\begin{eqnarray}
&& R^{\trp}\left(R {\bf z}+{\bf a}+{\bf d}-\boldsymbol \lambda\right)=0\label{cp1}\\
&& R{\bf z}+{\bf a}\geq 0\label{cp2}\\
&& \boldsymbol\lambda\cdot (R{\bf z}+{\bf a})=0\label{cp3}\\
&& \boldsymbol\lambda\geq 0\label{cp4}.
\end{eqnarray}
See \cite[Chapter 12]{NW}.




Before stating and proving the main result (Theorem \ref{main-theorem}), we list all 
the assumptions, including both new and previously used.

\noindent
A1. Consider the problem of  minimizing (\ref{q-basic}) subject only to the boundary conditions
(\ref{b1}, \ref{b1-bis}), but not the constraints (\ref{i3-intro}). We assume that the minimizer
of that problem is not feasible, that is, this minimizer does not
satisfy the impenetrability constraints (\ref{i3-intro}). Further we assume that the feasible region is not
empty, which means that there is at lest one point ${\bf z}$ satisfying all the inequality constraints
(\ref{cp2}).

\noindent
A2. The network $\Gamma$ is a regular triangulation (as defined in Definition \ref{reg-treg}).

\noindent
A3. The boundary conditions are prescribed so that 
\begin{equation}
\label{upper}
|({\bf x}^i+{\bf u}^i)-({\bf x}^i+{\bf u}^j)|\leq a_i+a_j+\min_{k=1,2,
\ldots N}a_k,
\end{equation}
for each pair ${\bf x}^i$, ${\bf x}^j$ of {\it boundary vertices in contact}.

Let us provide some comments on the nature  of assumptions A1--A3.
Assumption A1 means that minimizing the energy of the spring network subject only to boundary conditions leads to a configuration in which at least one spring is compressed (and thus violates the impenetrability constraints).

Assumption A2 concerns the contact geometry. The edges of the network split
the domain of the problem (a polygon) into elementary cells (triangles). Near the boundary, the cells must be compatible with the geometry of the boundary in the following sense. If an interior vertex is connected to a boundary vertex, then it is also connected with another
boundary vertex, located next to the first boundary vertex. Hence, every triangular cell adjacent to the exterior boundary must contain one free vertex and two boundary vertices.

Assumption A3 means that the boundary conditions (\ref{b1}, \ref{b1-bis}) are chosen to prevent particles from escaping through the gaps made by displacing the boundary particles. Clearly, if two boundary particles belong to the same group, then no gap can appear between them, and $|({\bf x}^i+{\bf u}^i)-({\bf x}^i+{\bf u}^j)|=a_i+a_j$. Formation of gaps would be possible between two boundary
particles from different groups which are in contact in the reference configuration. If the parameters of rigid body motions in the boundary conditions (\ref{b1}, \ref{b1-bis}) are prescribed arbitrarily, then the two particles may move away from each other, and open a gap large enough for a third particle to slip through. Assumption A3 prohibits formation of such gaps.

\begin{theorem}
\label{main-theorem}
Suppose that assumptions A1-A3 hold. Then there exist $d^\star>0$ and a vector $\boldsymbol\delta=
(\delta_1, \delta_2, \ldots, \delta_E)\in {\reals}^E: \delta_l\in (-1, -1/2), l=1, 2, \ldots, E$, such that for each ${\bf d}=d\boldsymbol\delta$ with $d>d^\star$,
the unique minimizer of (\ref{pr1},\ref{pr2}) has the following property. Each
interior vertex ${\bf x}^i$ of $\Gamma$ has at least
two active edges incident at it. The corresponding unit vectors
${\bf q}^{i,j_1}, {\bf q}^{i,j_2}$ must be linearly independent.
\end{theorem}

\noindent
{\it Proof.}\newline
\noindent
{\it Step 1}. We claim that A3 implies that there is $c_0>0$, which
depends on the boundary conditions, but is independent
of the choice of $d_{ij}$ in (\ref{q-basic}), such that each feasible displacement
${\bf u}^i, i=1,2, \ldots, N$, satisfies
\begin{equation}
\label{up-u}
|u^i_k|\leq c_0,
\end{equation}
$k=1,2$. Indeed, first we observe that if
${\bf u}^i, {\bf u}^j$ satisfy the linearized constraint (\ref{i3-intro}),
then they also satisfy the distance constraint (\ref{0.1}) (the converse is not
true in general). Then any feasible collection of displacements also satisfies
the distance constraints (\ref{0.1}) for each pair of neighboring vertices. Now we recall the assumption made in the introduction to conclude that (\ref{0.1}) must hold for all pairs of vertices. Fix $l\in \{1,2,\ldots, N\}$, corresponding to an interior vertex, and consider a smaller packing ${\mathcal P}$ of particles, containing only $D_l$ and
all boundary particles. In the reference configuration, $D_l$ is completely surrounded by boundary particles. Then the boundary conditions are prescribed according
to A3, the boundary particles still completely confine $D_l$, so that
${\bf x}^l$ must displace to ${\bf x}^l+{\bf w}^l$ that lies inside
a certain bounded domain $\Omega^\prime$ that depends only on boundary conditions.
Since ${\bf x}^i_k$ are bounded, this implies that the claim is true for
all displacements ${\bf w}^l$ which are feasible for the smaller packing ${\mathcal P}$.
Clearly the set of all such displacements is larger than the set of all ${\bf u}^l$
feasible under all constraints (\ref{0.1}), and the latter set is larger than the set of all ${\bf u}^l$ feasible under the linearized constraints (\ref{i3-intro}).
This proves the claim.

\noindent 
{\it Step 2}.
Let 
\begin{equation}
\label{vi}
{\bf v}^i=\sum_{j\in{\mathcal N}_i}{\bf q}^{ij},~~~i=1,2,\ldots, 2(N-N_b).
\end{equation}
First, we prove the theorem under the additional assumption
\begin{equation}
\label{a4}
\mbox{For each}~i=1,2,\ldots, 2(N-N_b),~~{\bf v}^i\ne s{\bf q}^{ij},
\end{equation}
where
$j\in {\mathcal N}_i, s\in {\reals}$.

We note that (\ref{a4}) implies that 
\begin{equation}
\label{ii}
|{\bf v}^i|\geq v_0>0
\end{equation}
with $v_0$ independent of $i$. Indeed, $s$ can be zero, so validity of (\ref{a4}) means in particular
that all ${\bf v}^i$ are non-zero. Since there is finitely many ${\bf v}^i$, (\ref{ii}) holds.

Consider solutions of the KKT system
(\ref{cp1},\ref{cp2},\ref{cp3},\ref{cp4}). From (\ref{cp3}), (\ref{cp4})
it follows that $\lambda_j=0$ if the $j$-th constraint is inactive.
Let $\theta_j=(R{\bf z}+{\bf a})_j$.
If the $j$-th constraint is active then $\theta_j=0$,
while $\lambda_j$ is arbitrary. Suppose that a feasible point
${\bf z}^\star$ is given. Then
$\theta_j$ are given. To solve (\ref{cp1}) we need to find
$\boldsymbol \lambda$. Denote
by ${\bf r}^k, k=1,2,\ldots, E$
the rows of $R$ (the columns of $R\trp$), and suppose
that the rows ${\bf r}^1, {\bf r}^2\ldots, {\bf r}^S$ correspond to
the active constraints, and that the rows
${\bf r}^{S+1}, {\bf r}^{S+2},\ldots, {\bf r}^{E}$ correspond to the
inactive constraints. Choose ${\bf d}=(-d, -d, \ldots, -d)$. 
Then (\ref{cp1}) can be written as
\begin{equation}
\label{cp1-1}
-\sum_{l=1}^S {\bf r}^l \lambda_l+\sum_{l=S
+1}^E{\bf r}^l\theta_l+d\sum_{l=1}^E {\bf r}^l=0.
\end{equation}

Pick a vertex ${\bf x}^i$ of $\Gamma$ and consider the restriction
of each ${\bf r}^l$ in (\ref{cp1-1}) to the two components
corresponding to ${\bf x}^i$.
Then we have
\begin{equation}
\label{cp-local}
-\sum_{j\in {\mathcal N}_i}^\prime\lambda_{ij} {\bf q}^{ij}+\sum_{j\in {\mathcal N}_i}^{\prime\prime}\theta_{ij} {\bf q}^{ij}+
d{\bf v}^i=0,
\end{equation}
where the first sum is taken over active edges incident at ${\bf x}^i$, while the second sum is over the inactive edges
incident at ${\bf x}^i$. 

Next, we determine the minimal number of active edges needed for (\ref{cp-local}) to hold. We can look at
(\ref{cp-local}) as a local problem in which ${\bf u}^i$ may vary, while ${\bf u}^j, j\in {\mathcal N}_i$ are fixed.
Denote by ${\mathcal F}_i\subset \reals^2$ the feasible region of
this local problem.
By A3, ${\mathcal F}_i$ is a polygon, each side
of which corresponds to one or more constraints being active. 

In the generic case, one constraint per side is active. In the non-generic case, two or more
active constraints correspond to the same side. Since our goal is estimating the number of
active constraints from below, it is sufficient to consider only the generic case,
corresponding to the ``worst case scenario".
In the generic case there are only three possibilities.\newline
\noindent
{\it Case 1}. ${\bf u}^i$ is inside ${\mathcal F}_i$. All edges incident at ${\bf x}^i$ are inactive.\newline
\noindent
{\it Case 2}. ${\bf u}^i$ belongs to only one of the sides of $\partial{\mathcal F}_i$. One edge is active.\newline
\noindent
{\it Case 3}. ${\bf u}^i$ is a vertex of ${\mathcal F}_i$. Two edges are active.

Consider case 1.
Then (\ref{cp-local}) cannot hold for $d$ sufficiently large. Indeed, $|{\bf v}^i|\geq v_0>0$ by assumption, while
$|\sum_{j\in {\mathcal N}_i}\theta_{ij} {\bf q}^{ij}|$ is bounded from above independent of $d$ in view of (\ref{up-u}).

Consider case 2.
Let us number the active edge by $(i,1)$. Then (\ref{cp-local}) can be written as
\begin{equation}
\label{cp-loc-1}
-\lambda_{i1}{\bf q}^{i1}+\sum_{j\in {\mathcal N}_i, j>1} \left(({\bf u}^i-{\bf u}^j)\cdot{\bf q}^{ij}\right){\bf q}^{ij}+
d{\bf v}^i=0.
\end{equation}
Enlarging $d$, if necessary, we see that (\ref{cp-loc-1}) can hold only if
\begin{equation}
\label{crit1}  
{\bf v}^i=s{\bf q}^{i1},
\end{equation}
 where
$s<0$. Since
(\ref{crit1}) is not allowed by (\ref{vi}), (\ref{cp-local}) cannot hold for sufficiently large $d$.

Consider case 3. Number the two active edges by $(i, 1)$, $(i,2)$. The equation (\ref{cp-local}) is
\begin{equation}
\label{cp-loc-2}
-\lambda_{i1}{\bf q}^{i1}-\lambda_{i2}{\bf q}^{i2}+
\sum_{j\in {\mathcal N}_i, j>2} \left(({\bf u}^i-{\bf u}^j)\cdot{\bf q}^{ij}\right){\bf q}^{ij}+
d{\bf v}^i=0.
\end{equation}
For this to hold for large $d$, ${\bf v}^i$ must be a non-positive linear combination
of ${\bf q}^{i1}, {\bf q}^{i2}$. These two vectors are linearly independent, otherwise
their intersection would not be a vertex of ${\mathcal F}_i$. So, Case 3 is possible, provided
${\bf v}^i$ lies in the negative cone of two active edges.

\noindent
{\it Step 3}. Now we remove the assumption (\ref{a4}). For each $i=1, 2, \ldots, (N-N_b)$, and each $\boldsymbol\delta\in {\reals}^E$, define a two-dimensional
vector
$\tilde{\bf v}^i$ to be the restriction of $R^{\trp} \boldsymbol\delta\equiv\sum_{l=1}^E \delta_l{\bf r}^l$ to a place $i$. The theorem will be proved is we show that there is
a choice of $\boldsymbol\delta$ such that $\delta_{l}\in [1/2, 1]$, and $\tilde{\bf v}^i$ has property (\ref{a4}). 
Indeed, if such $\boldsymbol\delta$ is found, we could choose ${\bf d}=d\boldsymbol\delta$, where $d>0$ is sufficiently large, and repeat the arguments made in the first step, using 
$\tilde{\bf v}^i$ instead of ${\bf v}^i$.

To show existence of $\boldsymbol\delta$, consider the cube $C_E=\{{\bf y}\in {\reals}^E: y_l\in (1/2, 1), l=1, 2, \ldots, E\}$. Pick any point ${\bf y}^\star\in C_E$. Since $C_E$ is open, there is a Euclidean open ball
$B({\bf y}^\star, \rho)\subset C_E$, with the radius $\rho>0$. Consider the image of $B({\bf y}^\star, \rho)$
under the mapping $R^{\trp}$. Since $R$ has full rank, $R^{\trp}$ is surjective, and is therefore an open mapping.
Thus, $R^{\trp}(B({\bf y}^\star, \rho))$ contains a Euclidean open ball $B(R^{\trp}{\bf y}^\star, \rho^\star)$ of a positive radius $\rho^\star$ depending only on $R^{\trp}$ and $\rho$, but not on ${\bf y}^\star$. If $R^{\trp} {\bf y}^\star$ has
property (\ref{a4}), we choose $\boldsymbol\delta={\bf y}^\star$ and we are done. Otherwise,
note that for each $i=1, 2, \ldots, (N-N_b)$, the ball $B(R^{\trp}{\bf y}^\star, \rho^\star)$ contains a non-empty two-dimensional
Euclidean open ball $B_i$ centered at the restriction of $R^{\trp}{\bf y}^\star$ to the place $i$. Since for each $i$
the set $\{{\bf v}\in {\reals}^2: {\bf v}=s{\bf q}^{ij}, s\in {\reals}, j\in {\mathcal N}^i\}$ is a union of a finite number of lines, it cannot contain a two-dimensional ball. Therefore, for each $i=1, 2, \ldots, (N-N_b)$ 
there must be a vector $\tilde{\bf v}^i\in B_i$ having property (\ref{a4}). Now we can define 
$\sum_{l=1}^E \delta_l{\bf r}^l$ via its restrictions $\tilde{\bf v}^i$. Next, by construction, we can find
a vector $\boldsymbol\delta\in B({\bf y}^\star, \rho)\subset C_E$ such that $R^{\trp}\boldsymbol\delta=\sum_{l=1}^E \delta_l{\bf r}^l$. 

The theorem is proved.

\noindent
{\bf Remark}. 
The choice 
of $\boldsymbol\delta$ in the proof is based on the following criterion. Consider the vector $R^{\trp}\boldsymbol\delta\in {\reals}^{2(N-N_b)}$. The proof works if
\begin{equation}
\label{last}
R^{\trp}\boldsymbol\delta\ne({\bf v}_1, {\bf v}_2, \ldots, {\bf v}_{N-N_b}),
\end{equation}
where at least one of the two-dimensional vectors ${\bf v}_i, i=1, 2, \ldots, (N-N_b)$ is of the form
${\bf v}_i=s{\bf q}^{ij}$, for some real $s$ and $j\in {\mathcal N}_i$. In other words, if ${\bf v}_i$ 
is inadmissible, then it must lie on the line through the origin with direction vector ${\bf q}^{ij}$. 
For each fixed $i$, the inadmissible set $V_i=\{{\bf v}_i\in \reals^2: {\bf v}_i=s{\bf q}^{ij}\}$  has Hausdorff dimension one (a finite union of lines in the plane), while the admissible set is the two-dimensional complement of $V_i$ in $\reals^2$. Therefore, the set of inadmissible vectors  in the right hand side of (\ref{last}) has dimension $2(N-N_b)-1$ while the set of admissible vectors $R^{\trp}\boldsymbol\delta$ is of dimension $2(N-N_b)$. Thus the admissible $R^{\trp}\boldsymbol\delta$ are generic, and the theorem holds for a generic (in the sense of Definition \ref{reg-treg}) packing under a generic pre-stress (the latter is determined by a generic choice of $\boldsymbol\delta$).


\section{Order parameter}
Recently, a phenomenological theory of slow dense granular flows was proposed in
\cite{AT, VTA}. A key quantity in that theory is the order parameter, defined as the ratio
of the number of solid-like contacts to the number of all contacts within a given control volume.
In \cite{VTA}, a contact is considered solid-like if two particles are jammed together for longer
than a characteristic collision time. The relevant characteristic
time is $\tau=a/v_a$, where $a$ is particle radius and $v_a$ is the speed of sound in
a solid material of the particles. Our model corresponds to the instantaneous material
response, when $\tau$ is much smaller than other relaxation times in the system, such
as the ratio of the sample size to a typical particle velocity.

An obvious type of pair motion leading to a solid-like contact is a rigid displacement
(a pair of particles infinitesimally moves as a rigid body). We shall call this type of contacts {\it stuck}. If a contact between $D_i$ and $D_j$ is stuck, then $({\bf u}^i-{\bf u}^j)\cdot{\bf q}^{ij}=0$, which is easy to check using the definition of rigid
displacements. This means that the impenetrability constraint for the corresponding
edge of $(i,j)$ of the network is satisfied as an equation (the edge is active). 
However, not every
active edge corresponds to a stuck contact. Another type of a local motion that produces
$({\bf u}^i-{\bf u}^j)\cdot{\bf q}^{ij}=0$ is an infinitesimal shear motion when
${\bf u}^i-{\bf u}^j$ is orthogonal to ${\bf q}^{ij}$. The corresponding contact is
called {\it sheared}. Note also that infinitesimal shear
is the same as infinitesimal rotation, so this type of motion includes infinitesimal rolling
as well as shear sliding.

We consider both sheared and stuck contact as solid-like, because stuck contacts are stable, while sheared contacts in an actual granular material will be subject to friction. Friction can be viewed as partially stabilizing, at least when the shearing force
is below the static friction threshold. Such non-sliding frictional contacts are considered as solid-like in the simulations performed in \cite{VTA}. In addition, some heuristic arguments and numerical simulations 
presented in \cite{GG1, GG2}, suggest that friction enhances elastic behavior of sufficiently large samples. Therefore, it makes sense to think of the network of solid-like contacts as the
main load-bearing structure and call this network {\it strong}. 
In contrast, a broken contact satisfying (\ref{broken}) corresponds to
a local weakening in the material because in this case two particles separate completely.
We can think of the network of all broken contacts as {\it weak}. 
Moreover, division of contacts into broken and solid-like corresponds to
the division of constraints into active and inactive, as done in optimization theory.
Therefore, this division is natural mathematically, and also makes sense from the physics point of view.

In addition, the definition in \cite{VTA} does not sufficiently clarify the nature of averaging. 
The notion of an order parameter in static problems should not use time averaging. 
The result of spatial averaging depends on the size of the sample that is being averaged. Thus, if the order
parameter is obtained by, say, spatial averaging, then it must depend on both location and size of the 
``control volume". In the discrete situation, the size of the averaging sample can be measured
by the minimal number of edges connecting a pair of vertices within the sample.

This suggests a definition of the size-dependent order parameter. To state
this definition we first define the averaging sample. 
\begin{definition}
\label{k-n}
A vertex ${\bf x}^j$ is in the k-th neighborhood of ${\bf x}^i$
if $\Gamma$ contains a path connecting ${\bf x}^i$ and ${\bf x}^j$ with
no more than $k$ edges.
\end{definition}
Now, to each $k$-neighborhood we can associate a value of an order parameter.  
\begin{definition}
\label{order}
For each ${\bf x}^i$ and each non-negative integer $k\leq N$, the size dependent
order parameter $\rho({\bf x}^i, k)$ is defined by
\begin{equation}
\label{op}
\rho({\bf x}^i, k)=\frac{\sum_{k} n_s}{\sum_k n},
\end{equation}
where the numerator is the number of active edges in $k$-neighborhood of ${\bf x}^i$,
and denominator is the number of all edges in that neighborhood.
\end{definition}

Theorem \ref{main-theorem} implies the lower bound 
\begin{equation}
\label{simple}
\rho({\bf x}^i,N)\geq \frac{N}{E}
\end{equation}
on the order parameter associated with the maximal, $N$-th neighborhood of each interior vertex ${\bf x}^i$.
Indeed, counting active edges (two per vertex) gives $2N$ edges, each counted at most twice.
In particular, (\ref{simple}) means that the order parameter $\rho({\bf x}^i, N)$
is bounded from below by the reciprocal of the mean coordination number of the network.

\section{Conclusions}

We have studied a network model of quasi-static deformation of dense pre-stressed
granular materials. 
In our model, the packing was represented by a network of linear elastic springs.
Each spring corresponds to a contact between two particles.
Geometric impenetrability constraints within the packing were modeled by the
linearized impenetrability constraints on the displacements of the vertices of the
network. The constraints have the form of linear inequalities, that can be satisfied either
as an equality (an active constraint), or as a strict inequality (inactive constraint).
Constraints are in one-to one correspondence with the interpaticle contacts. An active constraint
corresponds to a relatively stable solid-like contact. Inactive constraints represent the relatively weak
broken contacts.
The question addressed in the paper is to estimate the total number and distribution of the solid-like
contacts in the energy-minimizing configuration. We showed that each interior vertex of the network has at least two 
solid-like contacts corresponding to it. This result qualitatively
reproduces the micro-band structure obtained in \cite{Kuhn1, Kuhn2} by numerical simulations. We also discussed the connection between our result and a lower bound on the
order parameter \cite{AT, VTA}. In the paper, we proposed a definition of the
order parameter that is similar to the one introduced in \cite{VTA}, but differs from it in the interpretation of the so-called solid-like contacts. On the one hand, our definition appears to be in accord with a physical picture of granular statics, recently
proposed in \cite{GG1, GG2}. On the other hand, it is a naturally related to optimization theory.


\end{document}